\magnification=1100
\font\bigbf=cmbx10 scaled \magstep2
\font\medbf=cmbx10 scaled \magstep1 
\hfuzz=20 pt
\input psfig.sty

\font\medbf=cmbx10 at 13pt
\def\i{\item}
\def\n{\noindent}

\def\L{{\bf L}}

\def\ve{\varepsilon}
\def\n{\noindent}

\def\c{\centerline}
\def\vp{\varphi}

\def\vp{\varphi}

\def\O{{\cal O}}

\def\R{I\!\!R}

\def\vs{\vskip 2em}
\def\vsk{\vskip 4em}
\def\v{\vskip 1em}
\def\sqr#1#2{\vbox{\hrule height .#2pt
\hbox{\vrule width .#2pt height #1pt \kern #1pt
\vrule width .#2pt}\hrule height .#2pt }}
\def\square{\sqr74}
\def\endproof{\hphantom{MM}\hfill\llap{$\square$}\goodbreak}

\null
\vsk
\c{\bigbf On the Blow-up for a Discrete Boltzmann Equation}
\v
\c{\bigbf in the Plane}
\vsk
\c{Alberto Bressan and Massimo Fonte}
\v
\c{S.I.S.S.A., Via Beirut 4, Trieste 34014, ITALY.}
\vsk\vsk
\n{\bf Abstract.} We study the possibility of finite-time blow-up
for a two dimensional Broadwell model. In a set of rescaled variables,
we prove that no self-similar blow-up solution exists, and derive
some a priori bounds on the blow-up rate.  In the final section, 
a possible blow-up scenario is discussed.
\vsk\vsk

\

\n{\medbf 1 - Introduction}
\v
Consider the simplified model of a gas whose particles 
can have only finitely many speeds, say 
$c_1,\ldots,c_N\in\R^n$.  Call $u_i=u_i(t,x)$ the density
of particles with speed $c_i$.
The evolution of these densities can then be described by a semilinear
system of the form
$$\partial_t u_i +c_i\cdot \nabla u_i =\sum_{j,k} a_{ijk}\, u_ju_k
\qquad\qquad i=1,\ldots,N.
\eqno(1.1)$$
Here the coefficient $a_{ijk}$ measures the rate at which new
$i$-particles are created, as a result of collisions between
$j$- and $k$-particles. 
In a realistic model, these coefficients must satisfy a set of identities,
accounting for the conservation of mass, momentum and energy.

Given a continuous, bounded  initial data
$$u_i(0,x)=\bar u_i(x),\eqno(1.2)$$
on a small time interval $t\in [0,T]$
a solution of the Cauchy problem can be constructed by the
method of characteristics.
Indeed, since the system is semilinear,
this solution is obtained as the fixed point of the
integral transformation
$$u_i(t,x)=\bar u_i(x-c_it)+\int_0^t
\sum_{j,k} a_{ijk}\, u_ju_k \big(s,~x-c_i(t-s)\big)\,ds\,.\eqno(1.3)$$
For sufficiently small time intervals, the existence of a unique
fixed point follows from the contraction mapping principle, 
without any assumption on the constants $a_{ijk}$.

If the initial data is suitably small, the solution remains uniformly 
bounded for all times [3]. For large initial data, on the
other hand, the global existence and stability of solutions is 
known only in the one-dimensional case [2, 6, 10]. 
Since the right hand side has quadratic growth,
it might happen that the solution blows up in finite time.
Examples where the $\L^\infty$ norm of the solution
becomes arbitrarily large as $t\to\infty$ are easy to construct
[7].  In the present paper we focus on the
two-dimensional Broadwell model and examine the possibility that
blow-up actually occurs in finite time. 

Since the equations (1.1) admit a natural symmetry group,
one can perform an asymptotic rescaling of variables
and ask whether there is a blow-up solution which, in the rescaled variables,
converges to a steady state.   This technique has been widely used to
study blow-up singularities of reaction-diffusion equations with
superlinear forcing terms [4, 5].  See also [9] for an example of
self-similar blow-up for hyperbolic conservation laws.
Our results show, however, that for the two-dimensional Broadwell model
no such self-similar blow-up solution exists. 

If blow-up occurs at a time $T$, our results imply that for times $t
\to T-$ one has
$$\big\|u(t)\big\|_{\L^\infty}~>~ {1\over 5}\, {\ln \big|\ln(T-t)\big|
\over T-t}\,.\eqno(1.4)$$
This means that the blow-up rate must be different from
the natural growth rate $\big\|u(t)
\big\|_{\L^\infty} = \O(1)\cdot(T-t)^{-1}$ 
which would be obtained
in case of a quadratic equation $\dot u=C\,u^2$.

In the final section of this paper we discuss a possible scenario for
blow-up.  The analysis highlights how carefully chosen should
be the initial data, if blow-up is ever to happen.
This suggests that finite time blow-up is a highly non-generic
phenomenon, something one would not expect to encounter 
in  numerical simulations.
\vsk
\n{\medbf 2 - Coordinate rescaling}
\v
In the following, we say that $P^*=(t^*,x^*)$ is a {\bf blow-up point}
if 
$$\limsup_{~x\to x^*,\, t\to t^*-} u_i(t,x)=\infty$$
for some $i\in\{1,\ldots,N\}$.  
Define the constant
$$C\doteq \max_i |c_i|\,.$$
We say that
$(t^*,x^*)$ is a {\bf primary blow-up point}
if it is a blow-up point and the backward cone
$$\Gamma\doteq \big\{ (t,x)\,;~~|x-x^*|< 2C\,(t^*-t)\big\}$$
does not contain any other blow-up point.
\v
\n{\bf Lemma 1.} {\it Let $u=u(t,x)$ be a solution
of the Cauchy problem (1.1)-(1.2) with continuous initial data.
If no primary blow-up point exist, then $u$ is continuous
on the whole domain $[0,\infty[\,\times\R^n$.}
\v
\n{\bf Proof.} If $u$ is not continuous, it must be unbounded in
the neighborhood of some point. Hence some blow-up point 
exists.
Call ${\cal B}$ the set of such blow-up points. Define the function
$$\vp(x)\doteq \inf_{(\tau,\xi)\in{\cal B}}\big\{\tau+C\,|x-\xi|\big\}\,.
$$
By Ekeland's variational principle (see [1], p.254), 
there exists a point $x^*$ such that
$$\vp(x)\geq \vp(x^*)-{C\over 2}\, |x-x^*|$$
for all $x\in\R^2$.  Then $P^*\doteq \big(\vp(x^*),\,x^*\big)$
is a primary blow-up point. 
\endproof
Let now $(t^*,x^*)$ be a primary blow-up point.
One way to study the local asymptotic behavior of
$u$ is to rewrite the system in
terms of the rescaled variables $w_i= w_i(\tau,\eta)$, defined by
$$\left\{\eqalign{\tau &~=~ -\ln (t^*-t),\cr
\eta&~=~e^{\tau}x~=~{x-x^*\over {t^*-t}}\,,
\cr
w_i &~=~e^{-\tau}u_i~=~(t^*-t)u_i.\cr}\right.\eqno(2.1)$$
The corresponding system of evolution equations is
$$\partial_\tau w_i+(c_i+\eta)\cdot \nabla_\eta w_i
=~-w_i+\sum_{j,k} a_{ijk} \, w_jw_k
\qquad\qquad i=1,\ldots,n.\eqno(2.2)$$
Any nontrivial stationary or periodic solution $w$ of (2.2)
would yield a solution $u$ of (1.1) which blows up at $(t^*,x^*)$.
On the other hand, the non-existence of such solutions
for (2.2) would suggest that finite time blow-up
for (1.1) is unlikely.
\vsk

\n{\medbf 3 - The two-dimensional Broadwell model}
\v
Consider a system on $\R^2$
consisting of 4 types particles (fig.~1), with speeds
$$c_1=(1,1),\qquad c_2=(1,-1),\qquad c_3=(-1,-1),\qquad c_4=(-1,1).$$

\midinsert
\vskip 10pt
\centerline{\hbox{\psfig{figure=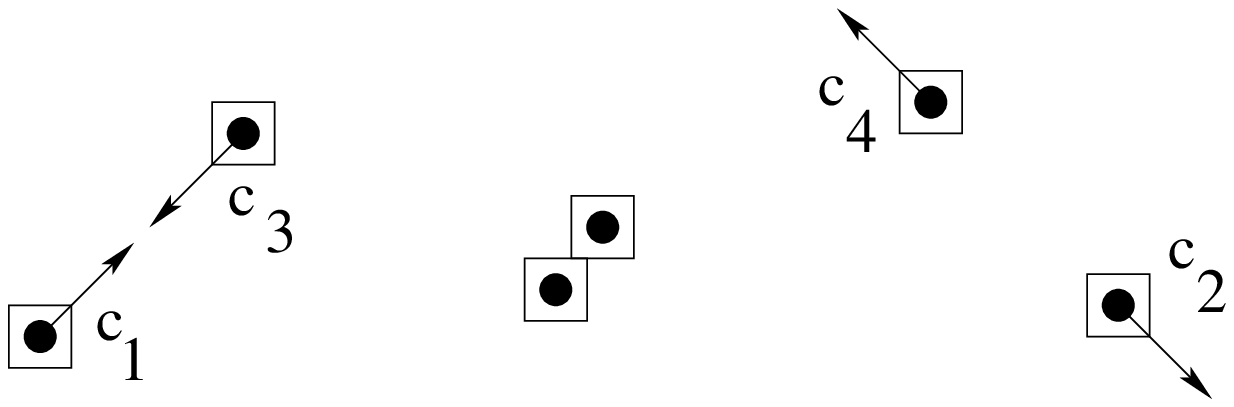,width=10cm}}}
\centerline{figure 1}
\vskip 10pt
\endinsert

The evolution equations are
$$\left\{\eqalign{\partial_t u_1+c_1\cdot\nabla u_1=u_2u_4-u_1u_3\,,\cr
\partial_t u_3+c_3\cdot\nabla u_3=u_2u_4-u_1u_3\,,\cr
\partial_t u_2+c_2\cdot\nabla u_2=u_1u_3-u_2u_4\,,\cr
\partial_t u_4+c_4\cdot\nabla u_4=u_1u_3-u_2u_4\,.\cr}\right.\eqno(3.1)$$
After renaming variables,
the corresponding rescaled system (2.2) takes the form
$$\left\{\eqalign{
\partial_t w_1+(x+1)\partial_x w_1+(y+1)\partial_y w_1
=w_2w_4-w_1w_3-w_1\,,\cr
\partial_t w_3+(x-1)\partial_x w_3+(y-1)\partial_y w_3
=w_2w_4-w_1w_3-w_3\,,\cr
\partial_t w_2+(x+1)\partial_x w_2+(y-1)\partial_y w_2
=w_1w_3-w_2w_4-w_2\,,\cr
\partial_t w_4+(x-1)\partial_x w_4+(y+1)\partial_y w_4
=w_1w_3-w_2w_4-w_4\,.\cr}
\right.\eqno(3.2)$$
\v
Our first result rules out the possibility of asymptotically
self-similar blow-up solutions.  A sharper estimate will be proved later.
\v
\n{\bf Theorem 1.} {\it The system (3.2) admits no 
nontrivial positive bounded solution which is constant or periodic in time.
}
\v
\n{\bf Proof.}  Assume 
$$0\leq w_i(t,x,y)\leq \kappa\eqno(3.3)
$$
for all  $t,x,y$, $i=1,2$.
Choose $\ve\doteq e^{-2\kappa}/2$, so that
$$\ve<{1\over\kappa}\,,\qquad\qquad \ve e^{2\kappa x}\leq {1\over 2}
\quad x\in [-1,1]\,.$$
Define
$$Q_{14}(t,y)\doteq
\int_{-1}^1 \Big[\big(1-\ve e^{2\kappa x}\big) w_1(t,x,y)
+\big(1-\ve e^{-2\kappa x}\big) w_4(t,x,y)\Big]\,dx\,.$$
$$Q_{14}(t)\doteq \sup_{|y|\leq 1} Q_{14}(t,y)\,,$$

\midinsert
\vskip 10pt
\centerline{\hbox{\psfig{figure=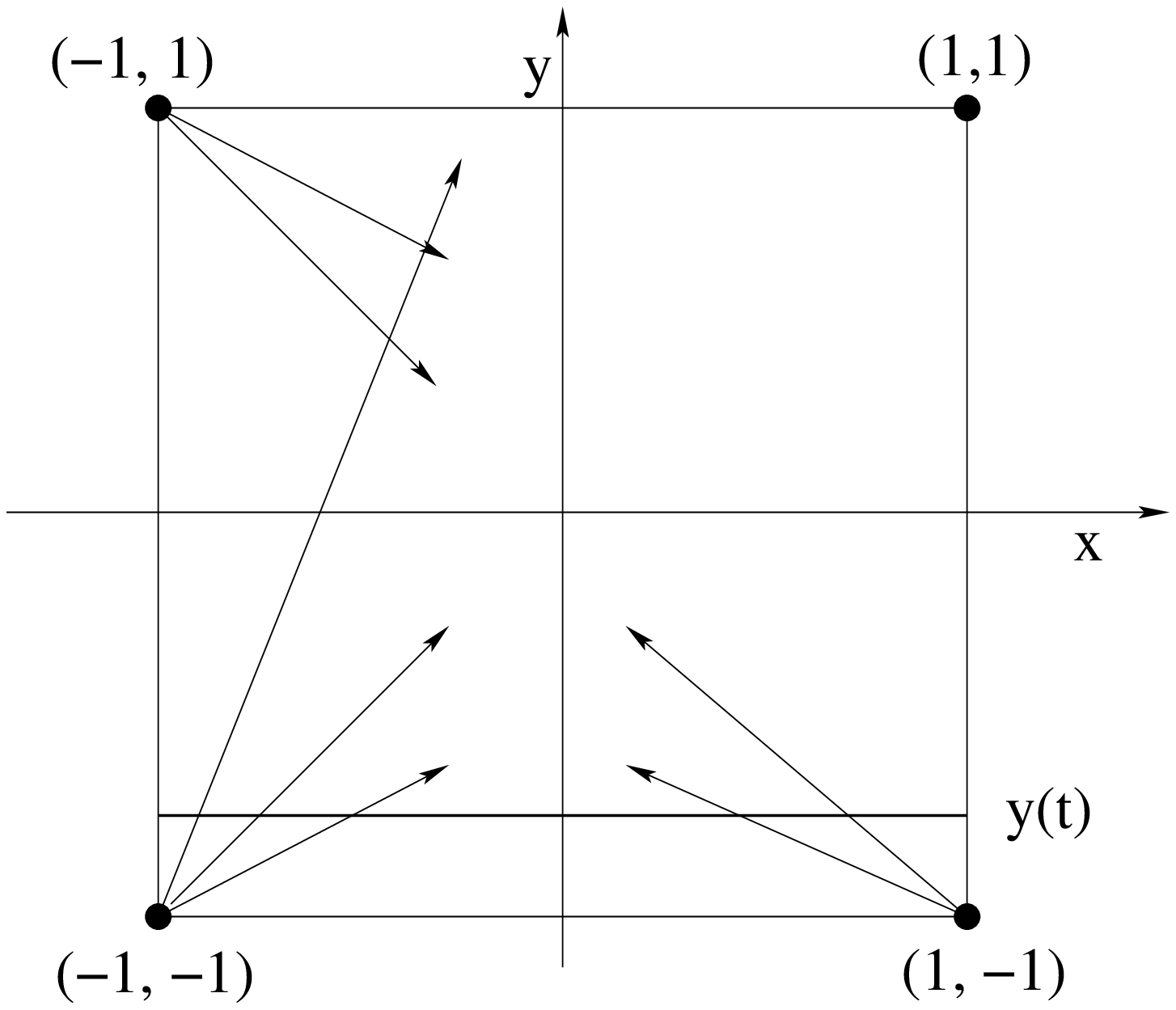,width=10cm}}}
\centerline{figure 2}
\vskip 10pt
\endinsert

Restricted to any horizontal moving line 
$y=y(t)$ such that $\dot y=y+1$ ~(fig.~2), the equations (3.2) become
$$\eqalign{
\partial_t w_1+(x+1)\partial_x w_1
=w_2w_4-w_1w_3-w_1\,,\cr
\partial_t w_4+(x-1)\partial_x w_4
=w_1w_3-w_2w_4-w_4\,.\cr}
$$
A direct computation now yields
$$\eqalign{{d\over dt}&Q_{14}\big(t,y(t)\big)\cr
&\leq
-2\ve\kappa \int_{-1}^1 \Big[e^{2\kappa x}(1+x)w_1
+e^{-2\kappa x}(1-x) w_4\Big] \,dx
+\int_{-1}^1 \big(\ve
e^{2\kappa x}-\ve
e^{-2\kappa x}\big)
(w_1w_3-w_2w_4)\,dx\cr
&\leq -\ve\kappa  \int_{-1}^1 \Big[e^{2\kappa x}(1+x)w_1
+e^{-2\kappa x}(1-x) w_4\Big] \,dx
-\ve\kappa \int_{-1}^0 e^{-2\kappa x}(1-x) w_4\,dx\cr
&\qquad 
-\ve\kappa \int_0^1e^{2\kappa x}(1+x)w_1\,dx
+\int_{-1}^0 \ve \kappa e^{-2\kappa x} w_4\,dx+
\int_0^1 \ve \kappa e^{2\kappa x} w_1\,dx\cr
&\leq -\ve\kappa  \int_{-1}^1 \Big[e^{2\kappa x}(1+x)w_1
+e^{-2\kappa x}(1-x) w_4\Big] \,dx\,.
}$$
Call 
$$I(t,y)\doteq\int_{-1}^1 w_1(t,x,y)\,dx\,.$$
The definition of $\ve$ and the bound (3.3) on $w_1$ imply 
$$\eqalign{\int_{-1}^1 e^{2\kappa x}(1+x)w_1\,dx
&\geq 2\ve \int_{-1}^1 (1+x)w_1\,dx\cr
&\geq  2\ve\int_{-1}^{-1+I/\kappa} (1+x)\kappa\,dx\cr
&=\ve\,I^2/\kappa\cr}$$
From this, and a similar estimate for $w_4$, we obtain
$$\int_{-1}^1 \Big[e^{2\kappa x}(1+x)w_1
+e^{-2\kappa x}(1-x) w_4\Big] \,dx\geq
{\ve\over\kappa}\left(\int_{-1}^1 w_1\,dx\right)^2+
{\ve\over\kappa}\left(\int_{-1}^1 w_4\,dx\right)^2\geq {\ve\over\kappa}\,
{Q^2_{14}\over 2}\,.$$
Since $\ve<\kappa^{-1}$, this yields
$$
{d\over dt}Q_{14}\big(t,y(t)\big)
\leq -{\ve^2\over 2}\,Q^2_{14}\big(t,y(t)\big)
\,.\eqno(3.4)$$
Observing that the Cauchy problem
$$\dot z=-{\ve^2\over 2}\,z^2,\qquad\qquad z(0)=4\kappa$$
has the solution  
$$z(t)=\left( {1\over 4\kappa}+{\ve^2\over 2}t\right)^{-1},$$
by a comparison argument from (3.4) we deduce
$$Q_{14}(t)\leq\left( {1\over 4\kappa}+{\ve^2\over 2}t\right)^{-1}.$$
Since
$$\int_{-1}^1\!\int_{-1}^1 w_1(t,x,y)\,dxdy\leq 4Q_{14}(t),$$
and since a similar estimate can be performed for all components
$w_i$, we conclude
$$\int_{-1}^1\!\int_{-1}^1 w_i(t,x,y)\,dxdy\leq 4
\left( {1\over 4\kappa}+{e^{-4\kappa}\over 8}t\right)^{-1}.
\eqno(3.5)
$$
The right hand side of (3.5) approaches zero as $t\to\infty$.
Therefore, nontrivial constant or time-periodic $\L^\infty$
solutions of (3.2) cannot exist.
\endproof
\vsk
\n{\medbf 4 - Refined blow-up estimates}
\v
If $(t^*,x^*)$ is a blow-up point, our analysis has 
shown that
in the rescaled coordinates $\tau,\xi$ the corresponding functions
$w_i$ must become unbounded as $\tau\to\infty$.   In this section
we refine the previous result, establishing a lower bound
for the rate at which such explosion takes place.
\v
\n{\bf Theorem 2.} {\it Let $u$ be a continuous
solution of the Broadwell system (2.1). Fix any point
$(t^*,x^*)$ and consider the corresponding rescaled variables
$\tau,\xi,w_i$. If
$$\max_{|\xi_1|, |\xi_2|\leq 1}\,
w_i(\tau,\xi_1,\xi_2)\leq \theta\ln\tau\qquad\qquad i=1,2,3,4\,,\eqno(4.1)$$
for some $\theta<1/4$ and all $\tau$ sufficiently large,
then
$$\lim_{\tau\to\infty} w_i(\tau,\xi)=0\qquad\qquad i=1,2,3,4\,,\eqno(4.2)$$
uniformly for $\xi \in\R^2$ in compact sets.
Therefore $(t^*,x^*)$ is not a blow up point.}
\v
Since $w_i=(t^*-t) u_i$ and 
$\tau\doteq \big|\ln(t^*-t)\big|$, the above implies
\v
\n{\bf Corollary 1.} {\it If $(t^*,x^*)$ is a primary blow-up point, 
then}
$$\limsup_{x\to x^*,~~t\to t^*- }~~  \big|u(t,x)\big|\cdot {
t^*-t\over  \ln\big|\ln(t^*-t)\big|}~\geq ~{1\over 4}\,.$$
\v
\n{\bf Proof of Theorem 2.} 
\v
Let $w_i=w_i(t,x,y)$ provide a solution to the system (3.2),
with 
$$0\leq w_i(t,x,y)\leq \theta \,\ln t\doteq k(t)\eqno(4.3)$$
for all $t\geq t_0$ and $x,y\in [-1,1]$. The proof will
be given in two steps. First we show that the $\L^1$ norm of the
components $w_i$ approaches zero as $t\to\infty$.
Then we refine the estimates, and prove that also the $\L^\infty$
norm asymptotically vanishes.
\v
\n STEP 1: Integral estimates.
Consider the function
$$
Q_{14}(t,y)\doteq\int_{-1}^1 \left[\left(1- {e^{2k(t)(x-1)}\over 2}\right)
w_1(t,x,y)+\left(1- {e^{-2k(t)(x+1)}\over 2}\right)w_4(t,x,y)\right]dx
$$
with $k(t)$ as in (4.3).
As in the proof of Theorem 1, let $t\mapsto y(t)$ be a solution to
$\dot y=y+1$. Then
$$
\eqalign{ {d\over{dt}}Q_{14}\big(t,\,y(t)\big)
=&\int_{-1}^1\big[-(x-1) k'  e^{2k(t)(x-1)}w_1+(x+1) k' e^{-2k(t)(x+1)}w_4
\big]\,dx
\cr
&+\int_{-1}^1\left(1-{{e^{2k(t)(x-1)}}\over 2}\right)
\big[-(x+1)w_{1x}+w_2w_4-w_1w_3-w_1\big]\,dx\cr
&+\int_{-1}^1\left(1-{{e^{-2k(t)(x+1)}}\over 2}\right)\big[-(x-1)w_{4x}+
w_1w_3-w_2w_4-w_4\big]\,dx}
$$
To estimate the right hand side, we notice that
$$
\eqalign{
A&\doteq\int_{-1}^1\left(1-{{e^{2k(t)(x-1)}}\over 2}\right)
\big[(1+x)w_{1x}+w_1\big]dx\geq k(t)\int_{-1}^1(x+1)e^{2k(x-1)}w_1\,dx
\cr
B&\doteq\int_{-1}^1\left(1-{{e^{-2k(t)(x+1)}}\over 2}\right)
\big[(x-1)w_{4x}+w_4\big]dx\geq k(t)\int_{-1}^1(1-x)e^{-2k(x+1)}w_4\,dx
\cr
C&\doteq\int_{-1}^1(w_1w_3-w_2w_4)\left({{e^{2k(x-1)}}\over 2}-
{{e^{-2k(x+1)}}\over 2}\right)\,dx
\cr
&\leq k(t)\int_0^1{{e^{2k(x-1)}}\over 2}w_1dx+k(t)\int_{-1}^0
{{e^{-2k(x+1)}} \over 2}w_4\,dx\,.}
$$
Therefore,
$$
\eqalign{{d\over{dt}}Q_{14}\big(t,\,y(t)\big)
=&\int_{-1}^1\big[-(x-1) k'e^{2k(t)(x-1)}w_1+(x+1) k' e^{-2k(t)(x+1)}
w_4\big]\,dx
-A-B+C\cr\leq&~k'(t)\int_{-1}^1\big[(1-x)e^{2k(t)(x-1)}w_1+(1+x)
e^{-2k(t)(x+1)}
w_4\big]\,dx
\cr
&\quad -{{k(t)}\over 2}\int_{-1}^1 \big[(1+x) e^{2k(t)(x-1)}w_1+(1-x)
e^{-2k(t)(x+1)}
w_4\big]\,dx.}$$
If $k(t)\geq 1/2$,
we claim that the following two inequalities hold: 
$$
\eqalign{
(1-x)e^{2k(t)(x-1)}\leq 1-e^{2k(t)(x-1)}\,,
\cr
(1+x)e^{-2k(t)(x+1)}\leq 1-e^{-2k(t)(x+1)}\,.
}\eqno(4.4)
$$
To prove the first inequality we need to show that
$$h_k(s)\doteq 
1-e^{2ks}+se^{2ks}\geq 0\qquad\qquad \hbox{for all}~~s\in [-2,0]\,.$$
This is clear because
$
h_k(0)=0
$
and
$$
h'_k(s)=e^{2ks}(1-2k+2ks)\leq 0\qquad\qquad s\in [-2,0]
$$
if $k\geq {1/ 2}$. Hence 
$h_k(s)$ is positive for $s\in [-2,0]$, as claimed.
The second inequality in (4.4) is proved similarly.

When $t\geq t_0\doteq e^{1/(2\theta)}$ one has
$k(t)\geq {1\over 2}$ and hence
$${d\over{dt}}Q_{14}\big(t,\,y(t)\big)
\leq
k'(t)Q_{14}- {{k(t)}\over 2}\int_{-1}^1 [(1+x) e^{2k(t)(x-1)}w_1+(1-x)
e^{-2k(t)(x+1)}w_4]dx\,.
$$
Setting $I=\int_{-1}^1w_1dx$, we 
obtain
$$
\int_{-1}^1(1+x)e^{2k(t)(x-1)}w_1dx 
~\geq~\int_{-1}^{-1+I/k(t)}(1+x)e^{-4k(t)} k(t)\,dx 
~=~e^{-4k(t)}{I^2\over 
{2k(t)}}\,.
$$
Using the above, and a similar estimate for the integral of $w_4$, we obtain
$$
\eqalign{
{{k(t)}\over 2}&\int_{-1}^1 [(1+x) e^{2k(t)(x-1)}w_1+(1-x)e^{-2k(t)(x+1)}w_4]
dx
\cr
&\geq{{ e^{-4k(t)}}\over 4}\left[\left( \int_{-1}^1w_1dx\right)^2+
\left( \int_{-1}^1w_4dx\right)^2\right]
\cr
&\geq {{ e^{-4k(t)}}\over 8}\,Q^2_{14}\,.\cr }
\eqno(4.5)$$
Calling 
$$Q_{14}(t)=\max_{|y|\leq 1} \,Q_{14}(t,y)\,,$$
from (4.5) we deduce
$$
{d\over dt}Q_{14}(t)
\leq
k'(t)Q_{14}(t)-{{ e^{-4k(t)}}\over 8}\,Q_{14}(t)^2.
$$
Recalling that $k(t)=\theta\ln t$ for some $0<\theta<1/4$, the previous 
differential inequality can be written as
$$
{d\over dt} Q_{14}
\leq  {\theta \over t} Q_{14}-{1\over {8t^{4\theta}}}Q_{14}^2.
\eqno(4.6)
$$
Notice that $ Q_{14}\big(t_0,y(t_0)\big)\leq 2k(t_0)$, and define
the constant
$$
A_0\doteq  \max\big\{2k(t_0)t_0^{1-4\theta}, ~8(1-3\theta)\big\}\,.$$
Then the function
$$
z(t)\doteq A_0\,t^{4\theta-1}
$$
satisfies
$$
{d\over dt} z(t)\geq {\theta \over t} z-{1\over {t^{4\theta}}}z^2\qquad 
\qquad z(t_0)\geq Q_{14}
(t_0,y(t_0))\,.\eqno(4.7)
$$
Comparing  (4.6) with (4.7) we conclude
$$Q_{14}(t)\leq z(t)\qquad\qquad t\geq t_0\,.\eqno(4.8)$$
This implies the estimate  
$$
\int_{-1}^1 w_i(t,x,y_0)\,dx \leq 2 Q_{14}(t)\leq 2 A_0\,t^{4\theta-1}
$$
for $t\geq t_0$, $i\in \{1,4\}$ and any $y_0\in [-1,1]$. 
An entirely similar
argument applied to $Q_{12}$, $Q_{23},\ldots$ yields the estimates
$$
\int_{-1}^1 w_i(t,x,y_0)dx \leq 2 A_0t^{4\theta-1}\,,\qquad\qquad
\int_{-1}^1 w_i(t,x_0,y)dy \leq 2 A_0t^{4\theta-1}\,.\eqno(4.9)$$
for $i=1,2,3,4$, $x_0,y_0\in [-1,1]$ and $t\geq t_0$.
\vs

\n STEP 2: Pointwise estimates. 
Using the integral bounds (4.9), we now seek a uniform bound of the
form
$$w_i(t,x,y)\leq C_0\eqno(4.10)$$
for some constant $C_0$ and all $x,y\in [-1,\,1]$, $t>0$.  
\v
To prove (4.10), let $t\mapsto x(t)$, $t\mapsto y(t)\big)\in [-1,1]$ 
be solutions of
$$\dot x=x+1\,,\qquad\qquad \dot y=y+1\,.$$
Call 
$$A(t)\doteq \int_{x(t)}^1 (w_1+w_4)\big(t,x,y(t)\big)\,dx\,.$$
From our previous  estimates (4.9) it trivially follows 
$$
A(t)\leq 4 A_0t^{4\theta-1}.
\eqno(4.11)
$$
The time derivative of
$A(t)$ is computed as
$$
\eqalign{
{{dA}\over {dt}}=&-(x(t)+1)(w_1+w_4)(t,x(t),y(t))+\int_{x(t)}^1
\big[\partial_t w_1+(y(t)+1)\partial_y w_1+\partial_tw_4+(y(t)+1)
\partial_y 
w_4\big]\,dx
\cr
=&-\big(x(t)+1\big)(w_1+w_4)\big(t,x(t),y(t)\big)-\int_{x(t)}^1\big[w_1+(x+1)
\partial_x w_1+w_4+(x-1)\partial_xw_4\big]\,dx
\cr
=&-(x(t)+1)(w_1+w_4)\big(t,x(t),y(t)\big)-\int_{x(t)}^1(w_1+w_4)dx
\cr
&-2w_1\big(t,1,y(t)\big)+(x(t)+1)w_1\big(t,x(t),y(t)\big)+
\int_{x(t)}^1w_1\,dx
\cr
&-2w_4\big(t,1,y(t)\big)+(x(t)-1)w_4\big(t,x(t),y(t)\big)+
\int_{x(t)}^1w_4\,dx
\cr
\leq& \big[x(t)-1-(x(t)+1)\big]w_4\big(t,x(t),y(t)\big)~=~-2w_4\big(
t,x(t),y(t)\big)
\,.\cr}
$$
This implies 
$$w_4\big(t,x(t),y(t)\big)\leq -{1\over 2} {{dA}\over{dt}}\,.\eqno(4.12)$$
The total derivative of $w_1$ along a characteristic line is now
given by
$$
\eqalign{{d\over {dt}}w_1\big(t,x(t),y(t)\big)
=&w_2w_4-w_1w_3-w_1\leq w_2w_4-w_1\leq {1\over 2}w_2\left({{-dA}\over{dt}}
\right)-w_1
\cr
\leq&-w_1+{k(t)\over 2}\left({{-dA}\over{dt}}\right)\,.\cr}
$$
In turn,  for $t\geq t_0$ this yields the inequality
$$
\eqalign{
w_1\big(t,x(t),y(t)\big)
&\leq e^{-(t-t_0)}
\left[w_1(t_0)+\int_{t_0}^te^{s-t_0}k(s)(-A'(s))ds\right]
\cr
&\leq e^{-(t-t_0)} \big[w_1(t_0)+A(t_0)k(t_0)\big]
+e^{-(t-t_0)}
\int_{t_0}^t A(s)\big(e^{s-t_0}k(s)\big)'\, ds.}\eqno(4.13)
$$
The first term on the right hand side of (4.13)
approaches zero exponentially fast.
Concerning the second, we have
$$e^{-(t-t_0)}
\int_{t_0}^t A(s)\,e^{s-t_0}\big(k(s)+k'(s)\big)\,ds
\leq \int_{t_0}^t e^{-(t-s)}\,2A_0 s^{4\theta-1}\,\left(\theta
\ln s+
{\theta\over s}\right)\,ds\,.
$$
This also approaches zero as $t\to\infty$.
Repeating the same computations for all components, we conclude that
for some time $t_1$ sufficiently large there holds
$$w_i(t_1,x,y) < {1\over 2}\qquad\qquad \hbox{for all}~x,y\in [-1,1]\,.
\eqno(4.14)$$
By continuity, the inequalities in (4.14) remain valid for all
$x,y$ in a slightly larger square, say
$[-1-\epsilon,~1+\epsilon]$.
For $t\geq t_1$ we now define
$$M(t)\doteq \max\Big\{ w_i(t,x,y)\,;~~i=1,2,3,4,~~~x,y\in
[-1-e^{t-t_1}\epsilon\,,~1+e^{t-t_1}\epsilon]\Big\}\,.$$
{}From the equations (3.2) and (4.14) it now follows
$${d\over dt} M(t)\leq -M(t)+M^2(t)\leq {M(t)\over 2}\,,\qquad\qquad M(t_1)
\leq {1\over 2}\,.$$
$$M(\tau)\leq \left[ 1+e^{\tau-t_1}\left({1\over M_1}-1\right)\right]^{-1}
\leq e^{-\tau}\cdot e^{t_1}\qquad\qquad\hbox{for all}~~ \tau\geq t_1.$$
Returning to the original variables $u_i=e^\tau w_i$, this yields
$$u_i\leq e^{t_1}$$
in a whole neighborhood of the point $P^*=(t^*, x^*)$.  Hence 
$P^*$ is not a blow-up point.
\endproof
\vsk
\n{\medbf 5 - A tentative blow-up scenario}
\v
For a solution of the rescaled equation (3.1), the total mass
$$m(t)\doteq \int_{-1}^1\int_{-1}^1 \sum_{i=1}^4 w_i(t,x,y)\, dxd y$$
may well become unbounded as $t\to\infty$.
On the other hand, the one-dimensional 
integrals along horizontal or vertical segments 
decrease monotonically.  
Namely, if $t\mapsto y(t) $ satisfies $\dot y=y+1$, then
$${d\over dt}\int_{-1}^1 \big[w_1(t,x,y(t))+w_4(t,x,y(t))\big]\,dx \leq 
0\,.
$$
Similarly, if $\dot x=x-1$, then 
$${d\over dt}\int_{-1}^1 \big[w_3(t,x(t),y)+w_4(t,x(t),y)\big]\,dy 
\leq 0\,.$$
Analogous estimates hold for
the sums $w_1+w_2$ and $w_2+w_3$.
Therefore, a bound on the initial data
$$w_i(0,x,y)\in [0, M]\qquad\qquad\hbox{for all}~~x,y\in [-1,1]\,,$$
yields uniform integral bounds on the line integrals of all
components:
$$\int_{-1}^1 w_i(t,x,y)\,dx\leq 4M\,,\qquad\qquad
\int_{-1}^1 w_i(t,x,y)\,dy\leq 4M\,.\eqno(5.1)$$

\midinsert
\vskip 10pt
\centerline{\hbox{\psfig{figure=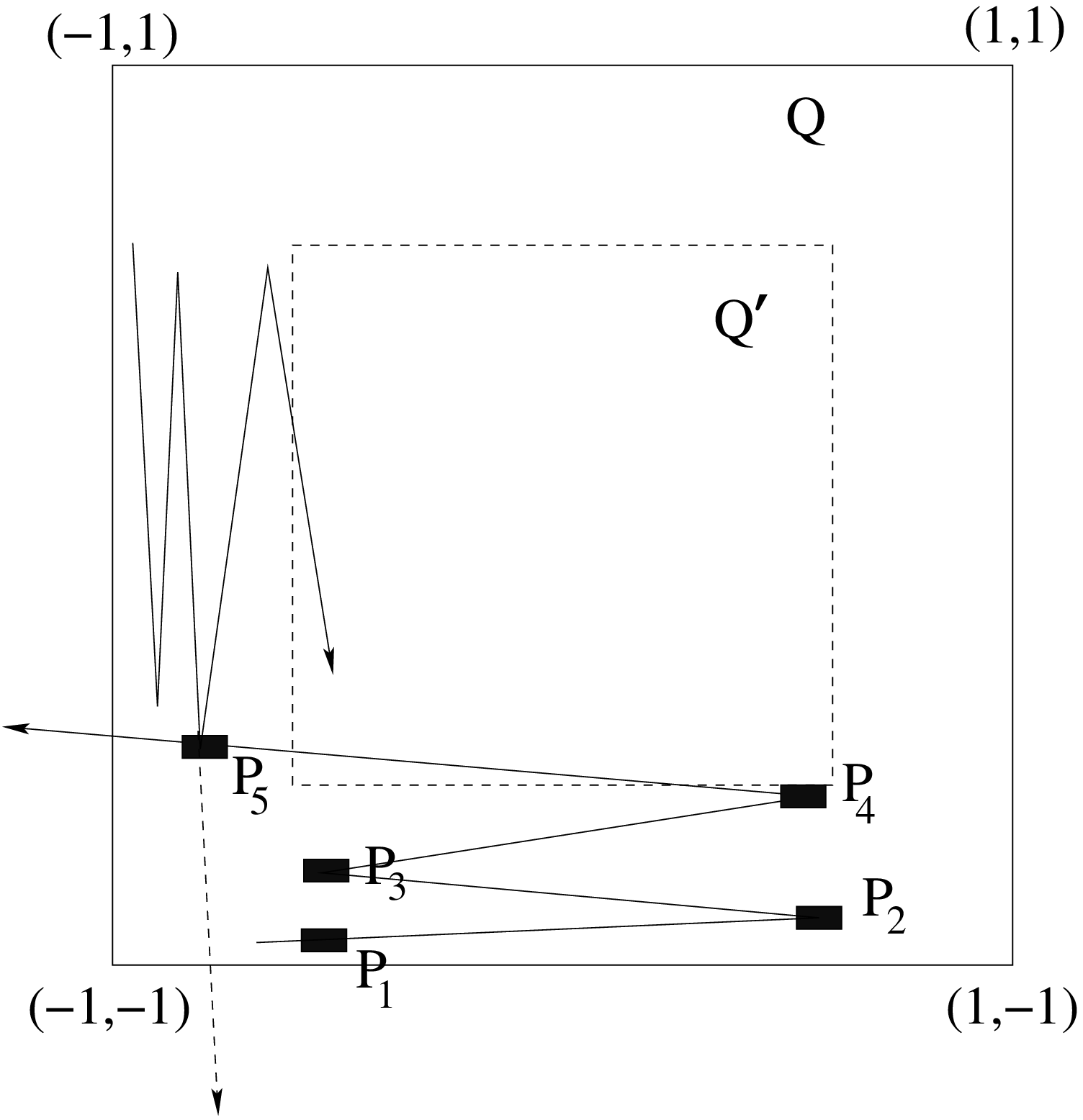,width=10cm}}}
\centerline{figure 3}
\vskip 10pt
\endinsert

If finite time blow-up is to occur, the mass which is initially 
distributed along each horizontal or vertical segment
must concentrate itself within a very small region, thus forming
a narrow packet of particles with increasingly high density.
A possible scenario is illustrated in fig.~3. 
A packet of 1-particles is initially located at $P_1$. 
In order to contribute to blow-up, this packet must remain
within the unit square $Q$. At $P_2$ these 1-particles interact
with 3-particles and produce a packet of 4-particles.
In turn, at $P_3$ these interact with 2-particles 
and produce again a packet of 1-particles.
After repeated interactions, the packet of alternatively 
1- and 4-particles eventually enters within the smaller square $Q'$.
After this time, it interacts with a packet of 2-particles at $P_5$
(transforming it into a packet of 1-particles)
and eventually exits from the domain $Q$.

To help intuition, it is convenient to describe a packet as
being ``young''
until it enters the smaller square $Q'$, and ``old'' afterwards.
To maintain a young packet inside $Q$, one needs the presence of
old packets interacting with it near the points
$P_2,P_3,P_4\ldots$ ~  On the other hand,
after it enters $Q'$, our packet can in turn be used to
hit another young packet, say at $P_5$, and preventing it from
leaving the domain $Q$.

As $t\to\infty$, the density of the packets must approach infinity.
One thus expects that most of the mass will be concentrated
along a finite number of one-dimensional curves.
Say, the packet of alternatively 1- and 4-particles should be
located along a moving curve $\gamma_{14}(t,\theta)$, where $\theta$
is a parameter along the curve.  The time evolution of such a curve
is of course governed by the equations

$${\partial\over\partial t}\gamma_{14}=c_1\qquad\hbox{or}\qquad
{\partial\over\partial t}\gamma_{14}=c_4$$
depending on whether $\gamma_{14}(t,\theta)$ consists of 1- or 4-particles.
The presence of interactions impose 
highly nonlinear constraints on
these curves.   For example,
the interaction occurring in $P_5$
at time $t$ implies
the crossing of the two curves $\gamma_{14}$ and $\gamma_{12}$,
namely
$$\gamma_{14}(t,\theta)=\gamma_{12}(t,\tilde \theta)=P_5$$
for some parameter values $\theta,\tilde\theta$.
The complicated geometry of these curves resulting from the
above constraints
has not been analyzed. 

\vsk
\c{\medbf References}
\v
\i{[1]} J.~P.~Aubin and I.~Ekeland, {\it Applied Nonlinear Analysis},
Wiley, 1984.
\v
\i{[2]} J.~M.~Bony, Solutions globales born\'ees pour les mod\`eles discrets
de l'\'equation de Boltzmann en dimension 1 d'espace,
{\it Actes Journ\'ees E.D.P. St. Jean de Monts} (1987).
\v
\i{[3]} J.~M.~Bony, Existence globale a donnee de Cauchy 
petites pour les modeles discrets de l'equation de Boltzmann, 
{\it Comm. Part. Diff. Equat.} {\bf 16} (1991), 533-545.
\v
\i{[4]} V.~A.~Galaktionov and J.~L.~Vazquez, The problem of blow-up
in nonlinear parabolic equations, {\it Discr. Cont. Dyn. Syst.} {\bf 8}
(2002), 399-433.
\v
\i{[5]} 
Y.~Giga and R.~Kohn, Characterizing blow-up using similarity variables,
{\it Indiana Univ. Math. J.} {\bf 36} (1987), 1-40.
\v
\i{[6]}
S.~Y.~Ha and A.~Tzavaras, Lyapunov functionals and
$L^1$ stability for discrete velocity Boltzmann equations, 
{\it Comm. Math. Phys.} {\bf 239} (2003), 65-92.
\v
\i{[7]} R.~Illner, Examples of non-bounded solutions in discrete 
kinetic theory, {\it J. M\'ecanique Th. Appl.} {\bf 5} (1986), 561-571.
\v
\i{[8]} R.~Illner and T. Platkowski, Discrete velocity models of the 
Boltzmann equation. A survey on the mathematical aspects of the theory,
{\it SIAM Review}{ \bf 30} (1988), 213-255.
\v
\i{[9]}  H.~K.~Jenssen,  Blowup for systems of conservation laws,
{\it SIAM J. Math. Anal.} {\bf 31} (2000), 894-908. 
\v
\i{[10]} 
L.~Tartar, Some existence theorem for semilinear hyperbolic 
systems in one space variable, {\it Technical Summary Report, Univ. 
Wisconsin} (1980).

\bye